\documentclass[11pt]{article}
\usepackage{latexsym,amsmath,amsfonts}
\begin{document}
\begin{center}
{\LARGE\textbf{Degree Distributions in General Random Intersection Graphs}}\\
\bigskip
\bigskip
Yilun Shang\footnote{Department of Mathematics, Shanghai Jiao Tong University, Shanghai 200240, CHINA. email: \texttt{shyl@sjtu.edu.cn}}\\
\end{center}

\begin{abstract}
We study a variant of the standard random intersection graph model
($G(n,m,F,H)$) in which random weights are assigned to both vertex
types in the bipartite structure. Under certain assumptions on the
distributions of these weights, the degree of a vertex is shown to
depend on the weight of that particular vertex and on the
distribution of the weights of the other vertex type.
\bigskip


\end{abstract}

\bigskip
\normalsize

\noindent{\Large\textbf{1. Introduction}}
\smallskip

Random intersection graphs, denoted by $G(n,m,p)$, are introduced in
\cite{9,3} as opposed to classical Erd\"os-R\'enyi random graphs.
Let us consider a set $V$ with $n$ vertices and another universal
set $W$ with $m$ elements. Define a bipartite graph $B(n,m,p)$ with
independent vertex sets $V$ and $W$. Edges between $v\in V$ and
$w\in W$ exist independently with probability $p$. The random
intersection graph $G(n,m,p)$ derived from $B(n,m,p)$ is defined on
the vertex set $V$ with vertices $v_1,v_2\in V$ adjacent if and only
if there exists some $w\in W$ such that both $v_1$ and $v_2$ are
adjacent to $w$ in $B(n,m,p)$.

To get an interesting graph structure and bounded average degree,
the work \cite{4} sets $m=\lfloor n^{\alpha}\rfloor$ and
$p=cn^{-(1+\alpha)/2}$ for some $\alpha$, $c>0$ and determines the
distribution of the degree of a typical vertex. Some related
properties for this model are recently investigated; for example,
independent set \cite{5} and component evolution \cite{7,6}. A
generalized random intersection graph is introduced in \cite{13} by
allowing a more general connection probability in the underlying
bipartite graph. The corresponding vertex degrees are also studied
by some authors, see e.g. \cite{12,10,11}, and shown to be
asymptotically Poisson distributed.

In this paper, we consider a variant model of random intersection
graphs, where each vertex and element are associated with a random
weight, in order to obtain tunable degree distributions. Our model,
referred to as $G(n,m,F,H)$, is defined as follows.
\smallskip

\noindent\textbf{Definition 1.}\itshape \quad Let us consider a set
$V=[n]$ of $n$ vertices and a set $W=[m]$ of $m$ elements. Define
$m=\lfloor\beta n^{\alpha}\rfloor$ with $\alpha,\beta>0$. Let
$\{A_i\}_{i=1}^n$ be an independent, identically distributed
sequence of positive random variables with distribution $F$. For
brevity, $F$ is assumed to have mean 1 if the mean is finite. The
sequence $\{B_i\}_{i=1}^m$ is defined analogously with distribution
$H$, which is independent with $F$ and assumed to have mean 1 if the
mean is finite. For some $i\in V$, $j\in W$ and $c>0$, set
\begin{equation}
p_{ij}=\big(cA_iB_jn^{-(1+\alpha)/2}\big)\wedge 1.\label{1}
\end{equation}
Define a bipartite graph $B(n,m,F,H)$ with independent vertex sets
$V$ and $W$. Edges between $i\in V$ and $j\in W$ exist independently
with probability $p_{ij}$. Then, $G(n,m,F,H)$ is constructed by
taking $V$ as the vertex set and drawing an edge between two
distinct vertices $i,j\in V$ if and only if they have a common
adjacent element $k\in W$ in $B(n,m,F,H)$.
 \normalfont\smallskip

If every element in $W$ has a unit weight, i.e. $H$ is a shifted
Heaviside function, our model reduces to that treated in \cite{1}.
Compared with Theorem 1.1 in \cite{1}, our result (see Theorem 1
below) provides more flexibility. A similar mechanism of assigning
random weights has been utilized for Erd\"os-R\'enyi graphs in
\cite{2} to generating random graphs with prescribed degree
distribution.

The rest of the paper is organized as follows. Our main results are
presented in Section 2 and we give proofs in Section 3.

\bigskip
\noindent{\Large\textbf{2. The results}}
\smallskip

Let $B$ be a random variable with distribution $H$ and suppose $B$
is independent with $\{B_i\}$. The following result concerns the
asymptotic expected degree of a vertex under appropriate moment
conditions on $F$ and $H$.

\smallskip
\noindent\textbf{Proposition 1.}\itshape\quad Let $D_i$ denote the
degree of vertex $i\in V$ in a general random intersection graph
$G(n,m,F,H)$ with $m=\lfloor\beta n^{\alpha}\rfloor$ and $p_{ij}$ as
in (\ref{1}). If $F$ has finite mean and $H$ has finite moment of
order 2, then, for all values of $\alpha>0$, we have that
$$
E(D_i|A_i)\rightarrow c^2A_i\beta E(B^2)
$$
almost surely, as $n\rightarrow\infty$. \normalfont
\smallskip

Our main theorem, which can be viewed as a generalization of Theorem
2 in \cite{4} and Theorem 1.1 in \cite{1}, reads as follows.

\smallskip
\noindent\textbf{Theorem 1.}\itshape\quad Let $D_i$ be the degree of
vertex $i\in V$ in a general random intersection graph $G(n,m,F,H)$
with $m=\lfloor\beta n^{\alpha}\rfloor$ and $p_{ij}$ as in
(\ref{1}). Assume that $F$ has finite mean.

(i) If $\alpha<1$, $H$ has finite moment of order
$(2\alpha/(1-\alpha))+\varepsilon$ for some $\varepsilon>0$, then,
as $n\rightarrow\infty$, the degree $D_i$ converges in distribution
to a point mass at 0.

(ii) If $\alpha=1$, $H$ has finite mean, then $D_i$ converges in
distribution to a sum of a $Poisson(cA_i\beta)$ distributed number
of $Poisson(cB)$ variables, where all variables are independent.

(iii) If $\alpha>1$, $H$ has finite moment of order 2, then $D_i$ is
asymptotically $Poisson(c^2A_i\beta B)$ distributed. \normalfont
\smallskip

The basic idea of proof is similar with that in \cite{1}, but some
significant modifications and new methods are adopted to tackle the
non-homogeneous connection probability involved here.

\bigskip
\noindent{\Large\textbf{3. Proofs}}
\smallskip

Let $|S|$ denote the cardinality of a set $S$. Suppose $\{x_n\}$ and
$\{y_n\}$ are sequences of real numbers with $y_n>0$ for all $n$, we
write $x_n\sim y_n$ if $\lim_{n\rightarrow\infty}x_n/y_n=1$; and if
$X$ and $Y$ are two random variables, we write $X\stackrel{d}{=}Y$
for equivalence in distribution. Without loss of generality, we
prove the results for vertex $i=1$.

\medskip
\noindent\textbf{Proof of Proposition 1}. We introduce cut-off
versions of the weight variables. For $i=2,\cdots, n$, let
$A'_i=A_i1_{[A_i\le n^{1/4}]}$ and $A''_i=A_i-A'_i$. Let $D'_1$ and
$D''_1$ be the degrees of vertex 1 when the weights
$\{A_i\}_{i\not=1}$ are replaced by $\{A'_i\}$ and $\{A''_i\}$,
respectively; that is, $D'_1$ is the number of neighbors of 1 with
weight less than or equal to $n^{1/4}$ and $D''_1$ is the number of
neighbors with weight larger than $n^{1/4}$. For $j\in W$, write
$p'_{1j}$ and $p''_{1j}$ for the analog of (\ref{1}) based on the
truncated weights.

For $i\in V$ and $i\not=1$, we observe that
$$
1-\prod_{j=1}^m(1-p_{1j}p''_{ij})\le\sum_{j=1}^mp_{1j}p''_{ij}\le
cA_1n^{-(1+\alpha)/2}\sum_{j=1}^mB_jp''_{ij}.
$$
Hence, we have
$$
E(D''_1|A_1)=\sum_{i=2}^nE\Big(1-\prod_{j=1}^m(1-p_{1j}p''_{ij})\Big)\le
c\beta
A_1n^{(\alpha-1)/2}\sum_{i=2}^n\bigg(\frac{\sum_{j=1}^mB_jEp''_{ij}}{m}\bigg).
$$
Since $F$ and $H$ have finite means, it follows that
$(\sum_{j=1}^mB_j)/m\rightarrow EB_1=1$ almost surely, by the strong
law of large numbers, and
$$
Ep''_{ij}\le
cn^{-(1+\alpha)/2}EA''_iEB_j=cn^{-(1+\alpha)/2}P(A_i>n^{1/4})\le
cn^{-(1+\alpha)/2}\frac{EA_i}{n^{1/4}},
$$
by using the Markov inequality. Therefore,
$E(D''_1|A_1)\rightarrow0$ almost surely, as $n\rightarrow\infty$.

As for $D'_1$, we observe that
$$
1-\prod_{j=1}^m(1-p_{1j}p'_{ij})=c^2A_1A'_i\Big(\sum_{j=1}^mB_j^2\Big)n^{-(1+\alpha)}+O\Big(A_1^2A^{\prime2}_i\Big(\sum_{k\not=l,k,l=1}^mB_k^2B_l^2\Big)n^{-2(1+\alpha)}\Big),
$$
and therefore,
\begin{eqnarray}
E(D'_1|A_1)&=&c^2A_1\beta
n^{-1}\Big(\sum_{i=2}^nEA'_i\Big)\bigg(\frac{\sum_{j=1}^mE(B_j^2)}{m}\bigg)\nonumber\\
&
&+n^{-2(1+\alpha)}O\Big(A_1^2E(A^{\prime2}_i)\Big(\sum_{k\not=l,k,l=1}^mE(B_k^2)E(B_l^2)\Big)\Big).\label{10}
\end{eqnarray}
The first term on the right-hand side of (\ref{10}) converges to
$c^2A_1\beta E(B^2)$ almost surely as $n\rightarrow\infty$ since
$\big(\sum_{j=1}^mE(B_j^2)\big)/m\rightarrow E(B^2)$ and
$EA'_i=EA_iP(A_i\le n^{1/4})\rightarrow EA_i=1$. The fact that
$A^{\prime2}_i\le n^{1/2}$ implies the second term on the right-hand
side of (\ref{10}) is  $O(n^{-2(1+\alpha)}n^{1/2}m^2)=o(1)$. The
proof is thus completed by noting that $D_1=D'_1+D''_1$. $\Box$

\medskip
\noindent\textbf{Proof of Theorem 1}. Let $N_1=\{j\in W|\ j\
\mathrm{is}\ \mathrm{adjaent}\ \mathrm{to}\ 1\in V\ \mathrm{in}\
B(n,m,F,H)\}$. Therefore, (i) follows if we prove that
$P(|N_1|=0)\rightarrow1$ as $n\rightarrow\infty$ for $\alpha<1$.
Conditional on $A_1$, $B_1,\cdots,B_m$, we have
\begin{equation}
P(|N_1|=0|\ A_1,B_1,\cdots,
B_m)=\prod_{k=1}^{m}(1-p_{1k})=1-O\Big(\sum_{k=1}^mp_{1k}\Big).\label{2}
\end{equation}
From (\ref{1}) we observe that
$$
\sum_{k=1}^mp_{1k}\le\sum_{k=1}^mcA_1B_kn^{-(1+\alpha)/2}\le m
\max_{k}\{B_k\}cA_1n^{-(1+\alpha)/2}=\beta
cA_1n^{(\alpha-1)/2}\max_{k}\{B_k\}.
$$
By the Markov inequality, for $\eta>0$
\begin{eqnarray*}
P\big(n^{(\alpha-1)/2}\max_{k}\{B_k\}>\eta\big)&\le&
mP\big(n^{(\alpha-1)/2}B_k>\eta\big)\\
&=&\beta
n^{\alpha}P\big(n^{-\alpha+\varepsilon(\alpha-1)/2}B_k^{(2\alpha/(1-\alpha))+\varepsilon}>\eta^{(2\alpha/(1-\alpha))+\varepsilon}\big)\\
&\le&\frac{\beta
E(B_k^{(2\alpha/(1-\alpha))+\varepsilon})}{\eta^{(2\alpha/(1-\alpha))+\varepsilon}n^{\varepsilon(1-\alpha)/2}}
\end{eqnarray*}
It then follows immediately from (\ref{2}) that $P(|N_1|=0|\
A_1,B_1,\cdots, B_m)\rightarrow1$ in probability, as
$n\rightarrow\infty$. Bounded convergence then gives that
$P(|N_1|=0)=EP(|N_1|=0|\ A_1, B_1,\cdots, B_m)\rightarrow1$, as
desired.

Next, to prove (ii) and (iii), we first note that
$ED''_1\rightarrow0$ as is proved in Proposition 1. The inequality
$P(D''_1>0)\le ED''_1$ implies that $D''_1$ converges to zero in
probability, and then it suffices to show that the generating
function of $D'_1$ converges to that of the claimed limiting
distribution. We condition on the variables $A_1$ and $B$, which are
assumed to be fixed in the sequel. For $i=2,\cdots,n$, let
$X'_i=\{j\in W|\ j\ \mathrm{is}\ \mathrm{adjacent}\ \mathrm{to}\
\mathrm{both}\ i\in V\ \mathrm{and}\ 1\in V\ \mathrm{in}\
B(n,m,F,H)\}$. Then by definition, we may write
$D'_1=\sum_{i=2}^n1_{[|X'_i|\ge1]}$. Conditional on $N_1,
A'_2,\cdots,A'_n,B_1,\cdots,B_m$, it is clear that $\{|X'_i|\}$ are
independent random variables and $X'_i\stackrel{d}{=}
\mathrm{Bernoulli}(p'_{ij_1})+\cdots+\mathrm{Bernoulli}(p'_{ij_{|N_1|}})$,
where the Bernoulli variables involved here are independent and we
assume $N_1=\{j_1,\cdots,j_{|N_1|}\}\subseteq W$. For $t\in[0,1]$,
the generating function of $D'_1$ can be expressed as
\begin{eqnarray*}
E\big(t^{D'_1}\big)&=&E\Big(\prod_{i=2}^nE\big(t^{1_{[|X'_i|\ge1]}}\big|\ N_1,A'_2,\cdots,A'_n,B_1,\cdots,B_m\big)\Big)\\
&=&E\Big(\prod_{i=2}^n\big(1+(t-1)P(|X'_i|\ge1|\
N_1,A'_2,\cdots,A'_n,B_1,\cdots,B_m)\big)\Big).
\end{eqnarray*}
Observe similarly as in Proposition 1 that
$$
P(|X'_i|\ge1|\
N_1,A'_2,\cdots,A'_n,B_1,\cdots,B_m)=1-\prod_{k=1}^{|N_1|}(1-p'_{ij_k})=\sum_{k=1}^{|N_1|}p'_{ij_k}+O\Big(\sum_{k\not=
l,k,l=1}^{|N_1|}p'_{ij_k}p'_{ij_l}\Big).
$$
Thereby, we have
\begin{eqnarray*}
\lefteqn{\prod_{i=2}^n\big(1+(t-1)P(|X'_i|\ge1|\
N_1,A'_2,\cdots,A'_n,B_1,\cdots,B_m)\big)}\\
& &=\exp\Big((t-1)\sum_{i=2}^n\sum_{k=1}^{|N_1|}p'_{ij_k}+O\Big(\sum_{i=2}^n\sum_{k,l=1}^{|N_1|}p'_{ij_k}p'_{ij_l}\Big)\Big)\\
& &=\exp\Big((t-1)\sum_{i=2}^n\sum_{k=1}^{|N_1|}p'_{ij_k}\Big)+R(n),
\end{eqnarray*}
where
$$
R(n):=\exp\Big((t-1)\sum_{i=2}^n\sum_{k=1}^{|N_1|}p'_{ij_k}\Big)\cdot\Big(\exp\Big(O\Big(\sum_{i=2}^n\sum_{k,l=1}^{|N_1|}p'_{ij_k}p'_{ij_l}\Big)\Big)-1\Big).
$$
Note that $E\big(t^{D'_1}\big)\in[0,1]$ and
$\exp\big((t-1)\sum_{i=2}^n\sum_{k=1}^{|N_1|}p'_{ij_k}\big)\in[0,1]$
since $t\in[0,1]$. Thus we have $R(n)\in[-1,1]$.

We then aim to prove the following three statements.

(a)
$E\big(\exp\big((t-1)\sum_{i=2}^n\sum_{k=1}^{|N_1|}p'_{ij_k}\big)\big)\rightarrow
e^{cA_1\beta(\exp(c(t-1)B)-1)}$, if $\alpha=1$;

(b)
$E\big(\exp\big((t-1)\sum_{i=2}^n\sum_{k=1}^{|N_1|}p'_{ij_k}\big)\big)\rightarrow
e^{c^2A_1\beta B(t-1)}$, if $\alpha>1$;

(c) $R(n)\rightarrow0$ in probability, if $\alpha\ge 1$.\\
The above limits in (a) and (b) are the generating functions for the
desired compound Poisson and Poisson distributions in (ii) and (iii)
of Theorem 1, respectively. By the bounded convergence theorem, (c)
yields $E(R(n))\rightarrow0$, which together with (a) and (b)
concludes the proof.

For $\alpha=1$, we have $|N_1|\stackrel{d}{=}
\mathrm{Bernoulli}(p_{11})+\cdots+\mathrm{Bernoulli}(p_{1m})$ and
all $m$ variables involved here are independent. By employing the
strong law of large numbers, we get
$$
\sum_{k=1}^mp_{1k}=cA_1\beta\frac{\sum_{j=1}^mB_j}{\beta
n}\rightarrow cA_1\beta\quad a.e.
$$
Then the Poisson paradigm (see e.g.\cite{14}) readily gives
$|N_1|\stackrel{d}{=}\mathrm{Poisson}(cA_1\beta)$. We have
\begin{eqnarray}
E\Big(\exp\Big((t-1)\sum_{i=2}^n\sum_{k=1}^{|N_1|}p'_{ij_k}\Big)\Big)&=&E\Big(E\Big(\exp\Big((t-1)\sum_{i=2}^n\sum_{k=1}^{|N_1|}p'_{ij_k}\Big)\Big|\
A'_2,\cdots,A'_n\Big)\Big)\nonumber\\
&=&E\Big(\sum_{s=0}^m\exp\Big((t-1)\sum_{i=2}^n\sum_{k=1}^{s}p'_{ik}\Big)\cdot
P(|N_1|=s)\Big).\label{3}
\end{eqnarray}
Since for any $k$ it follows that $EA'_i\rightarrow EA_i=1$ and
$\sum_{i=2}^np'_{ik}=cB_k(\sum_{i=2}^nA'_i)/n\rightarrow cB_k$
almost surely,
$$
\sum_{s=0}^m\exp\Big((t-1)\sum_{i=2}^n\sum_{k=1}^{s}p'_{ik}\Big)\cdot
P(|N_1|=s)\sim
\sum_{s=0}^m\exp\Big((t-1)c\sum_{k=1}^sB_k\Big)e^{-cA_1\beta}\frac{(cA_1\beta)^s}{s!}.
$$
For $t\in[0,1]$, define
$$
\tau=\tau(t):=Ee^{c(t-1)B}=\int e^{c(t-1)x}\mathrm{d}H(x),
$$
and then we obtain
\begin{eqnarray*}
E\Big(\sum_{s=0}^m\exp\Big((t-1)\sum_{i=2}^n\sum_{k=1}^{s}p'_{ik}\Big)\cdot
P(|N_1|=s)\Big)&\sim&E\Big(\sum_{s=0}^m\exp\Big((t-1)c\sum_{k=1}^sB_k\Big)e^{-cA_1\beta}\frac{(cA_1\beta)^s}{s!}\Big)\\
&=&\sum_{s=0}^m\Big(\prod_{k=1}^sE\big(e^{(t-1)cB_k}\big)\Big)e^{-cA_1\beta}\frac{(cA_1\beta)^s}{s!}\\
&=&e^{-cA_1\beta}\sum_{s=0}^m\frac{(\tau cA_1\beta)^s}{s!}\\
&\rightarrow&e^{cA_1\beta(\tau-1)}=e^{cA_1\beta(\exp(c(t-1)B)-1)}
\end{eqnarray*}
as $n\rightarrow\infty$, since $B$ is fixed. Combing this with
(\ref{3}) gives (a).

For $\alpha>1$, we also have $|N_1|\stackrel{d}{=}
\mathrm{Bernoulli}(p_{11})+\cdots+\mathrm{Bernoulli}(p_{1m})$ and
all $m$ variables involved here are independent. From strong law of
large numbers, it yields
\begin{equation}
\sum_{k=1}^mp_{1k}=cA_1\beta
n^{(\alpha-1)/2}\frac{\sum_{j=1}^mB_j}{\beta n^{\alpha}}\sim
cA_1\beta n^{(\alpha-1)/2}\quad a.e.\label{4}
\end{equation}
and for any $k$,
\begin{equation}
\sum_{i=2}^np'_{ik}=cB_kn^{(1-\alpha)/2}\frac{\sum_{i=2}^nA'_i}{n}\sim
cB_kn^{(1-\alpha)/2}\quad a.e.\label{5}
\end{equation}
Note that
\begin{equation}
\sum_{k=1}^mp_{1k}^2=\frac{\beta
c^2A_1^2}{n}\cdot\frac{\sum_{j=1}^mB_j^2}{\beta
n^{\alpha}}\rightarrow0\quad a.e.\label{6}
\end{equation}
as $n\rightarrow\infty$, since $H$ has finite moment of order 2. By
(\ref{4}), (\ref{6}) and a coupling argument of Poisson
approximation (see Section 2.2 \cite{8}), we obtain
$|N_1|\stackrel{d}{=} \mathrm{Poisson}(cA_1\beta n^{(\alpha-1)/2})$.
From (\ref{5}) we get
\begin{eqnarray*}
\lefteqn{\sum_{s=0}^m\exp\Big((t-1)\sum_{i=2}^n\sum_{k=1}^{s}p'_{ik}\Big)\cdot
P(|N_1|=s)}\\
& &\sim
\sum_{s=0}^m\exp\Big((t-1)cn^{(1-\alpha)/2}\sum_{k=1}^sB_k\Big)e^{-cA_1\beta
n^{(\alpha-1)/2}}\frac{(cA_1\beta n^{(\alpha-1)/2})^s}{s!}.
\end{eqnarray*}
For $t\in [0,1]$, define
$$
\omega=\omega(t,n):=Ee^{c(t-1)Bn^{(1-\alpha)/2}}
$$
and since (\ref{3}) still holds in this case, we similarly have
\begin{eqnarray*}
\lefteqn{E\Big(\sum_{s=0}^m\exp\Big((t-1)\sum_{i=2}^n\sum_{k=1}^{s}p'_{ik}\Big)\cdot
P(|N_1|=s)\Big)}\\
&\sim&\sum_{s=0}^m\Big(\prod_{k=1}^sE\big(e^{(t-1)cn^{(1-\alpha)/2}B_k}\big)\Big)e^{-cA_1\beta n^{(\alpha-1)/2}}\frac{(cA_1\beta n^{(\alpha-1)/2})^s}{s!}\\
&=&e^{-cA_1\beta n^{(\alpha-1)/2}}\sum_{s=0}^m\frac{(\omega cA_1\beta n^{(\alpha-1)/2})^s}{s!}\\
&\sim&e^{cA_1\beta n^{(\alpha-1)/2}(\omega-1)}\\
&=&\exp\Big(cA_1\beta
n^{(\alpha-1)/2}\Big(\frac{\exp\big(c(t-1)n^{(1-\alpha)/2}B\big)-1}{c(t-1)n^{(1-\alpha)/2}B}\Big)c(t-1)n^{(1-\alpha)/2}B\Big)\\
&\rightarrow&e^{c^2A_1\beta B(t-1)},
\end{eqnarray*}
as $n\rightarrow\infty$, where the last limit holds since
$c(t-1)n^{(1-\alpha)/2}B\rightarrow0$. Thus (b) is proved.

It remains to show (c). First note it suffices to show
\begin{equation}
\sum_{i=2}^n\sum_{k,l=1}^{|N_1|}p'_{ik}p'_{il}\rightarrow0\quad
\mathrm{in}\ \mathrm{probability}\label{7}
\end{equation}
as $n\rightarrow\infty$. Recalling that $A'_i\le n^{1/4}$, we have
for $\alpha\ge1$ that
$$
\sum_{k,l=1}^{|N_1|}\sum_{i=2}^np'_{ik}p'_{il}\le\sum_{k,l=1}^{|N_1|}c^2n^{-(1+\alpha)}B_kB_l\sum_{i=2}^nA^{\prime2}_i\le
c^2n^{(1/2)-\alpha}\Big(\sum_{k=1}^{|N_1|}B_k\Big)^2.
$$
For any $\eta>0$, we have
$$
P\Big(n^{(1/4)-\alpha/2}\sum_{k=1}^{|N_1|}B_k>\eta\Big)\le\frac{E\big(\sum_{k=1}^{|N_1|}B_k\big)}{\eta
n^{(\alpha/2)-1/4}}=\frac{(E|N_1|)(EB_1)}{\eta
n^{(\alpha/2)-1/4}}\le\frac{c\beta}{\eta n^{1/4}}
$$
by using the Markov inequality, the Wald equation (see e.g.
\cite{14}), $E|N_1|\le c\beta n^{(\alpha-1)/2}$ and $EB_1=1$,
proving the claim (\ref{7}) as it stands. $\Box$

\end{document}